\documentclass{article}
\usepackage{amsthm}
\usepackage{amssymb}
\usepackage{amsxtra}

\newtheorem{thm}{Theorem}[section]
\newtheorem{lemma}{Lemma}[section]

\title{Distribution of the zeros of the Riemann zeta function in
longer intervals}
\author{Tsz Ho Chan}

\begin{document}
\maketitle
\begin{abstract}
In this paper, we extend the result of Fujii on the second moment
of $S(t+h) - S(t)$ to longer range of $h$ under the Riemann
Hypothesis and an quantitative form of the Twin Prime Conjecture.
\end{abstract}
\section{Introduction}
Throughout this article, we shall assume the Riemann Hypothesis
{\bf RH} of the Riemann zeta function $\zeta(s)$. Let
$$S(t) = \frac{1}{\pi} \mbox{arg} \zeta\Bigl(\frac{1}{2} + it
\Bigr),$$
\begin{equation*}
F(x,T) = \sum_{0 < \gamma, \gamma' \leq T}
         x^{i(\gamma - \gamma')} w(\gamma - \gamma')
\mbox{ with } w(u) = \frac{4}{4 + u^2}.
\end{equation*}
Here, $\gamma$ and $\gamma'$ run over the imaginary parts of the
non-trivial zeros of $\zeta(s)$. In [\ref{F}], Fujii applied
Goldston's result [\ref{G}] to obtain, under {\bf RH}
\begin{equation*}
\begin{split}
\int_{0}^{T} (S(t+h) - S(t))^2\, dt &= \frac{T}{\pi^2} \Bigl[
\int_{0}^{h\log{(T/2\pi)}} \frac{1-\cos{\alpha}}{\alpha} d\alpha
\\
+& \int_{1}^{\infty} \frac{F(\alpha)}{\alpha^2} \Bigl(1 -
\cos{(\alpha h \log{\frac{T}{2\pi}})}\Bigr) d\alpha \Bigr] + o(T)
\end{split}
\end{equation*}
where
$$F(\alpha) = F(\alpha,T) := \Bigl(\frac{T}{2\pi}
\log{\frac{T}{2\pi}}\Bigr)^{-1} F\Bigl((\frac{T}{2\pi})^\alpha,
T\Bigr),$$ and $0 < h = o(1)$.

To extend the above result to a longer range $0< h = O(1)$, one
needs to improve the error term. The first step in this direction
was accomplished in the author's [\ref{C1}] which improves the
error term of the second moment of $S(t)$ to $O(T/\log^2{T})$
under an quantitative form of the Twin Prime Conjecture {\bf TPC}
(see next section). We shall use the more precise estimates in
[\ref{C1}] and [\ref{C2}] to prove
\begin{thm}
\label{theorem1} Assume {\bf RH} and {\bf TPC}. Fix a large number
$A$.
\begin{equation*}
\begin{split}
& \int_{0}^{T} (S(t+h) - S(t))^2\, dt \\
=& \frac{T}{\pi^2} \int_{0}^{hL} \frac{1 - \cos{\alpha}}{\alpha}
d\alpha + \frac{T}{\pi^2} \int_{1}^{\infty} \frac{F(\alpha) -
F_h(\alpha)}{\alpha^2} d\alpha
\\
+& \frac{T}{\pi^2} \Bigl[\log{\log{2}} + C_0 - \sum_{m=2}^{\infty}
\sum_{p} \frac{1}{mp^m} \Bigr] (1 - \cos{(h\log{2})}) -
\frac{T}{\pi^2} \int_{0}^{h\log{2}} \frac{1 -
\cos{\alpha}}{\alpha} d\alpha \\
-& \frac{Th}{\pi^2} \int_{2}^{\infty} r(u)
\frac{\sin{(h\log{u})}}{u} du + \frac{T}{\pi^2}
\sum_{m=2}^{\infty} \sum_{p} \frac{1 - \cos{(hm\log{p})}}{m^2 p^m}
\\
+& \frac{Th \sin{(hL)}}{\pi^2 (4+h^2)L} + \frac{T}{\pi^2 L^2}
\frac{h^2(20+3h^2)}{4(4+h^2)^2} -\frac{3T}{2\pi^2 L^2}
\int_{1}^{\infty} \frac{F(\alpha)-F_h(\alpha)}{\alpha^4} d\alpha +
O\Bigl(\frac{T}{L^2}\Bigr)
\end{split}
\end{equation*}
for $0 < h \leq A$. The implicit constant in the error term may
depend on $A$. $F_h(\alpha)$ and $r(u)$ are given by (\ref{1.1})
and (\ref{3.1}) respectively. $C_0$ is Euler's constant.
\end{thm}

This prompts us to study
$$\int_{1}^{\infty} \frac{F(\alpha) -
F_h(\alpha)}{\alpha^2} d\alpha \mbox{ and } \int_{1}^{\infty}
\frac{F(\alpha)-F_h(\alpha)}{\alpha^4} d\alpha.$$ We have

\begin{thm}
\label{theorem2} Assume {\bf RH}. For $T$ sufficiently large,
$$0 \leq \int_{1}^{\infty} \frac{F(\alpha) - F_h(\alpha)}
{\alpha^2} d\alpha < 9,$$ and
$$0 \leq \int_{1}^{\infty} \frac{F(\alpha) - F_h(\alpha)}
{\alpha^4} d\alpha < 6.$$
\end{thm}

The author would like to thank Prof. Daniel Goldston for
suggesting this problem. Here and throughout this paper, $p$ will
denote a prime number. $\Lambda(n)$ is von Mangoldt's lambda
function. Also, we have $L = \log{\frac{T}{2\pi}}$.
\section{Preparations}
We shall use a strong quantitative form of the Twin Prime
Conjecture (abbreviated as {\bf TPC}): For any $\epsilon
> 0$,
$$\sum_{n=1}^{N} \Lambda(n) \Lambda(n+d) = {\mathfrak S}(d) N + O(N^{1/2+
\epsilon}) \mbox{ uniformly in } |d| \leq N.$$ ${\mathfrak S}(d) =
2\prod_{p>2}\bigl(1-\frac{1}{(p-1)^2}\bigr) \prod_{p|d, p>2}
\frac{p-1}{p-2}$ if $d$ is even, and ${\mathfrak S}(d) = 0$ if $d$
is odd.

\bigskip

Also, we need to generalize $F(x,T)$ to
\begin{equation}
\label{1.1} F_h(x,T) = \sum_{0 < \gamma, \gamma' \leq T}
\cos{\bigl((\gamma - \gamma' - h) \log{x}\bigr)} w(\gamma -
\gamma' - h).
\end{equation}
Note that $F_0(x,T) = F(x,T)$. From [\ref{C2}],
\begin{equation}
\label{2.1}
\begin{split}
F_h(x,T) =& \frac{T}{2 \pi} \Bigl[\frac{4 \cos{(h \log{x})}}{4 +
h^2} \log{x} - \frac{8h \sin{(h \log{x})}}{(4 + h^2)^2} \Bigr] \\
& + \frac{T}{2 \pi x^2} \Bigl[ \Bigl(\log{\frac{T}{2 \pi}}\Bigr)^2
- 2\log{\frac{T}{2 \pi}} \Bigr] + O(x \log{x}) + O\Bigl(
\frac{T}{x^{1/2 - \epsilon}} \Bigr)
\end{split}
\end{equation}
for $1 \leq x \leq \frac{T}{\log^3{T}}$ under {\bf RH}, and
\begin{equation}
\label{2.2} F_h(x,T) = \frac{T}{2\pi} \Bigl[\frac{4
\cos{(h\log{x})}}{4+h^2} \log{x} - \frac{8h
\sin{(h\log{x})}}{(4+h^2)^2} \Bigr] + O(x)
\end{equation}
for $\frac{T}{\log^3{T}} \leq x \leq T$ under {\bf RH} and {\bf
TPC}. Define
$$F_h(\alpha) = F_h(\alpha, T) := \Bigl(\frac{TL}{2\pi}\Bigr)^{-1}
F_h\Bigl((\frac{T}{2\pi})^\alpha, T\Bigr)$$ for $\alpha \geq 0$
and $F_h(-\alpha) = F_h(\alpha)$. So, (\ref{2.1}) and (\ref{2.2})
can be summarized as
\begin{equation}
\label{2.3} F_h(\alpha) = \left\{ \begin{array}{ll}
\frac{4\cos{(hL\alpha)}}{4+h^2} \alpha -
\frac{8h\sin{(hL\alpha)}}{(4+h^2)^2} L^{-1}  & \\
+ (\frac{T}{2\pi})^{-2\alpha} [L - 2] & \mbox{ if } 0 \leq
\alpha \leq 1 - \frac{3\log{\log{T}}}{\log{T}}, \\
+ O(\alpha T^{\alpha -1}) + O(T^{-(1/2-\epsilon)\alpha}L^{-1}), &
 \\
\frac{4\cos{(hL\alpha)}}{4+h^2} \alpha + O(\alpha T^{\alpha -1}
L^{-1}), & \mbox{ if } 1 - \frac{3\log{\log{T}}}{\log{T}} \leq
\alpha \leq 1. \end{array} \right.
\end{equation}

Recall from [\ref{G}],
$$k(u) = \left\{ \begin{array}{ll} \Bigl(\frac{1}{2u} -
\frac{\pi^2}{2} \cot{(\pi^2 u)}\Bigr)^2, & \mbox{ if } |u| \leq
\frac{1}{2\pi}, \\
\frac{1}{4u^2}, & \mbox{ if } |u| > \frac{1}{2\pi}.
\end{array} \right.$$ and $\hat{k}(u)$ denotes the Fourier transform of
$k(u)$. One can easily check that $k''(u)$ is bounded, piecewise
continuous and $\ll u^{-4}$ when $u > \frac{1}{2\pi}$. Also,
\begin{equation}
\label{2.4} k(0) = 0, k\Bigl(\frac{1}{2\pi}\Bigr) = \pi^2; k'(0) =
0,k'\Bigl({\frac{1}{2\pi}}^{-}\Bigr) = \pi^5 - 4\pi^3.
\end{equation}

We need the following lemmas.
\begin{lemma}
\label{lemma2.1} Let $x = (T/2\pi)^\beta$ with $\beta > 0$,
\begin{equation*}
\begin{split}
\sum_{0 < \gamma, \gamma' \leq T} \hat{k}((\gamma - \gamma')
\log{x}) =& \frac{T L}{4\pi^2 \beta} \int_{-\infty}^{\infty}
F(\alpha) k\Bigl(\frac{\alpha}{2\pi \beta}\Bigr) d\alpha +
\frac{\pi^2 T}{16 L} \frac{F(\beta)}{\beta^2} \\
&- \frac{T}{64 \pi^4 L \beta^3} \int_{-\infty}^{\infty} F(\alpha)
k''\Bigl(\frac{\alpha}{2\pi \beta}\Bigr) d\alpha.
\end{split}
\end{equation*}
\end{lemma}

Proof: This is essentially Lemma 2.6 of [\ref{C1}].

\begin{lemma}
\label{lemma2.2} Let $x = (T/2\pi)^\beta$ with $\beta > 0$,
\begin{equation*}
\begin{split}
\sum_{0 < \gamma, \gamma' \leq T} \hat{k}((\gamma - \gamma' - h)
\log{x}) =& \frac{T L}{4\pi^2 \beta} \int_{-\infty}^{\infty}
F_h(\alpha) k\Bigl(\frac{\alpha}{2\pi \beta}\Bigr) d\alpha +
\frac{\pi^2 T}{16 L} \frac{F_h(\beta)}{\beta^2} \\
&- \frac{T}{64 \pi^4 L \beta^3} \int_{-\infty}^{\infty}
F_h(\alpha) k''\Bigl(\frac{\alpha}{2\pi \beta}\Bigr) d\alpha.
\end{split}
\end{equation*}
\end{lemma}

Proof: This is just very similar to Lemma \ref{lemma2.1} above. We
use the fact that $k(u)$ is even.

\begin{lemma}
\label{lemma2.3}
\begin{equation*}
\begin{split}
\int_{0}^{\beta^-} \sin{(hL\alpha)} k''\Bigl(\frac{\alpha}{2\pi
\beta}\Bigr) d\alpha =& 2(\pi^6 - 4\pi^4) \beta \sin{(hL\beta)} -
4\pi^4 hL \beta^2 \cos{(hL\beta)} \\
&- (2\pi hL\beta)^2 \int_{0}^{\beta} \sin{(hL\alpha)}
k\Bigl(\frac{\alpha}{2\pi \beta}\Bigr) d\alpha.
\end{split}
\end{equation*}
\end{lemma}

Proof: Use integration by parts twice and (\ref{2.4}).

\begin{lemma}
\label{lemma2.4}
\begin{equation*}
\begin{split}
& \int_{0}^{\beta^-} \alpha \cos{(hL\alpha)}
k''\Bigl(\frac{\alpha}{2\pi \beta}\Bigr) d\alpha \\
=& 2 (\pi^6 - 6\pi^4) \beta^2 \cos{(hL\beta)} + 4 \pi^4 hL \beta^3
\sin{(hL \beta)} \\
&- 8\pi^2 hL\beta^2 \int_{0}^{\beta} \sin{(hL\alpha)}
k\Bigl(\frac{\alpha}{2\pi \beta}\Bigr) d\alpha - (2\pi hL\beta)^2
\int_{0}^{\beta} \cos{(hL\alpha)} k\Bigl(\frac{\alpha}{2\pi
\beta}\Bigr) d\alpha.
\end{split}
\end{equation*}
\end{lemma}

Proof: Use integration by parts twice and (\ref{2.4}) again.

\begin{lemma}
\label{lemma2.5}
$$\int_{0}^{\beta^-} \alpha k''\Bigl(\frac{\alpha}{2\pi \beta}
\Bigr) d\alpha = 2(\pi^6 - 6\pi^4) \beta^2.$$
\end{lemma}

Proof: Set $h=0$ in Lemma \ref{lemma2.4}.

\begin{lemma}
\label{lemma2.6}
\begin{equation*}
\begin{split}
\int_{\beta}^{1} \frac{\cos{(hL\alpha)}}{\alpha^3} d\alpha =&
\frac{\cos{(hL\beta)}}{2\beta^2} - \frac{\cos{(hL)}}{2} - \frac{hL
\sin{(hL\beta)}}{2 \beta} \\
&+ \frac{hL \sin{(hL)}}{2} - \frac{(hL)^2}{2} \int_{\beta}^{1}
\frac{\cos{(hL\alpha)}}{\alpha} d\alpha.
\end{split}
\end{equation*}
\end{lemma}

Proof: Use integration by parts twice.

\begin{lemma}
\label{lemma2.7}
\begin{equation*}
\begin{split}
\int_{\beta}^{1} \frac{\sin{(hL\alpha)}}{\alpha^4} d\alpha=&
\frac{\sin{(hL\beta)}}{3\beta^3} - \frac{\sin{(hL)}}{3} + \frac{hL
\cos{(hL\beta)}}{6 \beta^2} - \frac{hL \cos{(hL)}}{6} \\
-& \frac{(hL)^2 \sin{(hL\beta)}}{6\beta} + \frac{(hL)^2
\sin{(hL)}}{6} - \frac{(hL)^3}{6} \int_{\beta}^{1}
\frac{\cos{(hL\alpha)}}{\alpha} d\alpha.
\end{split}
\end{equation*}
\end{lemma}

Proof: Use integration by parts thrice.
\section{$S_4$ and $S_5$}

We shall follow Fujii [\ref{F}] closely. Let $x = (T/2\pi)^\beta$
with $0 < \beta < 1$. By Goldston's explicit formula of $S(t)$ in
[\ref{G}] under {\bf RH}, Fujii got (see p. 76 $\&$ 77 of
[\ref{F}])
$$\int_{0}^{T} (S(t+h) - S(t))^2\, dt = S_3 + S_4 + S_5 +
O\Bigl(x \frac{\log{\log{x}}}{\log{x}}\Bigr) + O(\log^3{T})$$
where
\begin{equation*}
\begin{split}
S_3 =& \frac{2}{\pi^2 \log{x}} \sum_{0< \gamma, \gamma' \leq T}
\bigl(\hat{k}((\gamma-\gamma')\log{x}) - \hat{k}((\gamma-\gamma' -
h)\log{x}) \bigr), \\
S_4 =& \frac{T}{\pi^2} \sum_{p \leq x} \frac{\cos{(h\log{p})}
-1}{p} \Bigl(f^2\bigl(\frac{\log{p}}{\log{x}}\bigr) -
2f\bigl(\frac{\log{p}}{\log{x}}\bigr) \Bigr), \\
S_5 =& \frac{T}{\pi^2} \sum_{m=2}^{\infty} \sum_{p^m \leq x}
\frac{\cos{(hm\log{p})} -1}{m^2 p^m}
\Bigl(f^2\bigl(\frac{m\log{p}}{\log{x}}\bigr) -
2f\bigl(\frac{m\log{p}}{\log{x}}\bigr) \Bigr).
\end{split}
\end{equation*}
Here $f(u) = \frac{\pi u}{2} \cot{(\frac{\pi u}{2})}$ and
$\hat{k}(u)$ is defined as in the previous section. We note that
with Euler's constant $C_0$,
\begin{equation}
\label{3.1}
\begin{split}
\sum_{p \leq u} \frac{1}{p} =& \log{\log{u}} + C_0 + \sum_{p}
\Bigl(\log{(1-\frac{1}{p})} + \frac{1}{p}\Bigr) + r(u), \\
r(u) \ll& \frac{\log{u}}{\sqrt{u}}
\end{split}
\end{equation}
under {\bf RH}. Consider
$$\Sigma_1 = \sum_{p \leq x} \frac{\cos{(h\log{p})}
-1}{p} f^2\Bigl(\frac{\log{p}}{\log{x}}\Bigr), \, \Sigma_2 =
\sum_{p \leq x} \frac{\cos{(h\log{p})} -1}{p}
f\Bigl(\frac{\log{p}}{\log{x}}\Bigr).$$ By partial summation,
\begin{equation*}
\begin{split}
\Sigma_1 =& \int_{2}^{x} \frac{\cos{(h\log{u})} - 1}{u\log{u}}
f^2\Bigl(\frac{\log{u}}{\log{x}}\Bigr) du + \int_{2^-}^{x}
(\cos{(h\log{u})} - 1) f^2\Bigl(\frac{\log{u}}{\log{x}}\Bigr)
dr(u) \\
=& \Sigma_{1,1} + \Sigma_{1,2}, \\
\Sigma_2 =& \int_{2}^{x} \frac{\cos{(h\log{u})} - 1}{u\log{u}}
f\Bigl(\frac{\log{u}}{\log{x}}\Bigr) du + \int_{2^-}^{x}
(\cos{(h\log{u})} - 1) f\Bigl(\frac{\log{u}}{\log{x}}\Bigr)
dr(u) \\
=& \Sigma_{2,1} + \Sigma_{2,2}.
\end{split}
\end{equation*}
Then
$$\Sigma_1 - 2\Sigma_2 = (\Sigma_{1,1} - 2\Sigma_{2,1}) +
(\Sigma_{1,2} - 2\Sigma_{2,2}).$$ As $f(0) = 1$, $f(1) = 0$ and
$f(u) = 1 + O(u^2)$, by integration by parts,
\begin{equation*}
\begin{split}
\Sigma_{1,2} - 2\Sigma_{2,2} =& -r(2^{-})(\cos{(h\log{2})} - 1)
\Bigl[f^2\bigl(\frac{\log{2}}{\log{x}}\bigr) -
2f\bigl(\frac{\log{2}}{\log{x}}\bigr) \Bigr] \\
&- \int_{2^-}^{x} r(u) \frac{\cos{(h\log{u})} - 1}{u \log{x}}
\Bigl[2f\bigl(\frac{\log{u}}{\log{x}}\bigr)
f'\bigl(\frac{\log{u}}{\log{x}}\bigr) -
2f'\bigl(\frac{\log{u}}{\log{x}}\bigr) \Bigr] du \\
&+ h \int_{2^-}^{x} r(u)
\Bigl[f^2\bigl(\frac{\log{u}}{\log{x}}\bigr) -
2f\bigl(\frac{\log{u}}{\log{x}}\bigr) \Bigr]
\frac{\sin{(h\log{u})}}{u} du \\
=& \Bigl[\log{\log{2}} + C_0 - \sum_{m=2}^{\infty} \sum_{p}
\frac{1}{mp^m} \Bigr] (1 - \cos{(h\log{2})}) \\
&- h \int_{2^-}^{\infty} r(u) \frac{\sin{(h\log{u})}}{u} du +
O\Bigl(\frac{1}{\beta^3 L^3}\Bigr)
\end{split}
\end{equation*}
by $f^2 - 2f = (f-1)^2 - 1$ and Taylor series of $\log{(1+x)}$.
\begin{equation*}
\begin{split}
& \Sigma_{1,1} - 2\Sigma_{2,1} \\
=& \int_{2}^{x} \frac{\cos{(h\log{u})} - 1}{u\log{u}} \Bigl[
f^2\bigl(\frac{\log{u}}{\log{x}}\bigr) -
2f\bigl(\frac{\log{u}}{\log{x}}\bigr) \Bigr] du\\
=& - \int_{2}^{x} \frac{1 - \cos{(h\log{u})}}{u\log{u}} \Bigl[1 -
\frac{\pi \log{u}}{2 \log{x}} \cot{\bigl(\frac{\pi \log{u}}{2
\log{x}}\bigr)} \Bigr]^2 du + \int_{2}^{x}
\frac{1 - \cos{(h\log{u})}}{u\log{u}} du \\
=& - \int_{0}^{\beta} \frac{1 - \cos{(hL\alpha)}}{\alpha} \Bigl[1
- \frac{\pi \alpha}{2 \beta} \cot{\bigl(\frac{\pi \alpha}{2 \beta}
\bigr)}\Bigr]^2 d\alpha + \int_{0}^{\beta} \frac{1 -
\cos{(hL\alpha)}}{\alpha} d\alpha \\
&- \int_{0}^{\log{2}/L} \frac{1 - \cos{(hL\alpha)}}{\alpha}
d\alpha + O\Bigl(\frac{1}{\beta^4 L^6}\Bigr)
\end{split}
\end{equation*}
by substituting $\alpha = \log{u}/L$. Therefore,
\begin{equation}
\label{S4}
\begin{split}
S_4 =& \frac{T}{\pi^2} \Bigl[\log{\log{2}} + C_0 -
\sum_{m=2}^{\infty} \sum_{p} \frac{1}{mp^m} \Bigr] (1 -
\cos{(h\log{2})}) \\
&+ \frac{T}{\pi^2} \int_{0}^{\beta} \frac{1 -
\cos{(hL\alpha)}}{\alpha} d\alpha - \frac{T}{\pi^2}
\int_{0}^{\beta} \frac{1 - \cos{(hL\alpha)}}{\alpha} \Bigl[1 -
\frac{\pi \alpha}{2 \beta} \cot{\bigl(\frac{\pi \alpha}{2 \beta}
\bigr)}\Bigr]^2 d\alpha \\
&- \frac{T}{\pi^2} \int_{0}^{h\log{2}} \frac{1 -
\cos{\alpha}}{\alpha} d\alpha - \frac{Th}{\pi^2}
\int_{2^-}^{\infty} r(u) \frac{\sin{(h\log{u})}}{u} du +
O\Bigl(\frac{1}{\beta^4 L^3}\Bigr).
\end{split}
\end{equation}
\begin{equation}
\label{S5}
\begin{split}
S_5 =& \frac{T}{\pi^2} \Bigl[\sum_{m=2}^{\infty} \sum_{p^m \leq x}
\frac{1 - \cos{(hm\log{p})}}{m^2 p^m} \\
&- \sum_{m=2}^{\infty} \sum_{p^m \leq x} \frac{1 -
\cos{(hm\log{p})}}{m^2 p^m}
\Bigl(f\bigl(\frac{m\log{p}}{\log{x}}\bigr) - 1\Bigr)^2 \Bigl]\\
=& \frac{T}{\pi^2} \sum_{m=2}^{\infty} \sum_{p} \frac{1 -
\cos{(hm\log{p})}}{m^2 p^m} + O\Bigl(\frac{1}{\beta^4 L^4}\Bigr).
\end{split}
\end{equation}
\section{$S_3$}

Applying Lemma \ref{lemma2.1} and Lemma \ref{lemma2.2}, we have
\begin{equation*}
\begin{split}
S_3 =& \frac{T}{2\pi^4 \beta^2} \int_{-\infty}^{\infty} (F(\alpha)
- F_h(\alpha)) k\Bigl(\frac{\alpha}{2\pi \beta}\Bigr) d\alpha +
\frac{T}{8 L^2} \frac{F(\beta) - F_h(\beta)}{\beta^3} \\
&- \frac{T}{32 \pi^6 L^2 \beta^4} \int_{-\infty}^{\infty}
(F(\alpha) - F_h(\alpha)) k''\Bigl(\frac{\alpha}{2\pi \beta}\Bigr)
d\alpha \\
=& S_{3,1} + S_{3,2} - S_{3,3}
\end{split}
\end{equation*}
From (\ref{2.3}),
\begin{equation}
\label{4.1} F(\alpha) - F_h(\alpha) = (1-\cos{(hL\alpha)})\alpha +
\frac{h^2}{4+h^2} \cos{(hL\alpha)} \alpha + \frac{8h
\sin{(hL\alpha)}}{(4+h^2)^2 L} + E
\end{equation}
where
\begin{equation}
\label{4.2} E = \left\{ \begin{array}{ll} O(T^{-(1/2 - \epsilon)
\alpha} L^{-1}) + O(\alpha T^{\alpha-1}), & \mbox{ if } 0 \leq
\alpha \leq
1 - \frac{3\log{\log{T}}}{\log{T}}, \\
O(T^{-(1/2 - \epsilon) \alpha} L^{-1}) + O(\alpha T^{\alpha-1}
L^{-1}), & \mbox{ if } 1 - \frac{3\log{\log{T}}}{\log{T}} \leq
\alpha \leq 1. \end{array} \right.
\end{equation}
Since $F(\alpha)$ and $F_h(\alpha)$ are even,
$$S_{3,1} = \frac{T}{\pi^4 \beta^2} \int_{0}^{\infty} (F(\alpha)
- F_h(\alpha)) k\Bigl(\frac{\alpha}{2\pi \beta}\Bigr) d\alpha.$$
Let $\epsilon_T = 3\log{\log{T}}/\log{T}$. We split the above
integral into four pieces:
$$I = \int_{0}^{\infty} = \int_{0}^{\beta} +
\int_{\beta}^{1-\epsilon_T} + \int_{1-\epsilon_T}^{1} +
\int_{1}^{\infty} = I_1 + I_2 + I_3 + I_4.$$ By (\ref{4.1}),
(\ref{4.2}) and the definition of $k(u)$,
\begin{equation*}
\begin{split}
I_1 =& (\pi \beta)^2 \int_{0}^{\beta} \frac{1 -
\cos{(hL\alpha)}}{\alpha} \Bigl[1 - \frac{\pi \alpha}{2 \beta}
\cot{\bigl(\frac{\pi \alpha}{2 \beta} \bigr)}\Bigr]^2 d\alpha \\
&+ \frac{h^2}{4+h^2} \int_{0}^{\beta} \alpha \cos{(hL\alpha)}
k\Bigl(\frac{\alpha}{2\pi \beta}\Bigr) d\alpha \\
&+ \frac{8h}{(4+h^2)^2 L} \int_{0}^{\beta} \sin{(hL\alpha)}
k\Bigl(\frac{\alpha}{2\pi \beta}\Bigr) d\alpha + O(\beta^{-1}
L^{-3}), \\
I_2 =& (\pi \beta)^2 \int_{\beta}^{1-\epsilon_T} \frac{1 -
\cos{(hL\alpha)}}{\alpha} d\alpha + \frac{(\pi \beta)^2 h^2}{4+h^2
} \int_{\beta}^{1-\epsilon_T} \frac{\cos{(hL
\alpha)}}{\alpha} d\alpha \\
&+ \frac{(\pi \beta)^2 8h}{(4+h^2)^2 L}
\int_{\beta}^{1-\epsilon_T} \frac{\sin{(hL \alpha)}}{\alpha}
d\alpha + O(L^{-4}), \\
I_3 =& (\pi \beta)^2 \int_{1-\epsilon_T}^{1} \frac{1 -
\cos{(hL\alpha)}}{\alpha} d\alpha + \frac{(\pi \beta)^2 h^2}{4+h^2
} \int_{1-\epsilon_T}^{1} \frac{\cos{(hL \alpha)}}{\alpha} d\alpha \\
&+ \frac{(\pi \beta)^2 8h}{(4+h^2)^2 L} \int_{1-\epsilon_T}^{1}
\frac{\sin{(hL \alpha)}}{\alpha} d\alpha + O(L^{-2}), \\
I_4 =& (\pi \beta)^2 \int_{1}^{\infty} \frac{F(\alpha) -
F_h(\alpha)}{\alpha^2} d\alpha.
\end{split}
\end{equation*}
Thus,
\begin{equation*}
\begin{split}
S_{3,1} =& \frac{T}{\pi^2} \int_{0}^{\beta} \frac{1 -
\cos{(hL\alpha)}}{\alpha} \Bigl[1 - \frac{\pi \alpha}{2 \beta}
\cot{\bigl(\frac{\pi \alpha}{2 \beta} \bigr)}\Bigr]^2 d\alpha +
\frac{T}{\pi^2} \int_{\beta}^{1} \frac{1 -
\cos{(hL\alpha)}}{\alpha} d\alpha \\
+& \frac{T}{\pi^4 \beta^2} \Bigl[\frac{h^2}{4+h^2}
\int_{0}^{\beta} \alpha \cos{(hL\alpha)} k\Bigl(\frac{\alpha}{2\pi
\beta}\Bigr) d\alpha \\
+& \frac{8h}{(4+h^2)^2 L} \int_{0}^{\beta} \sin{(hL\alpha)}
k\Bigl(\frac{\alpha}{2\pi \beta}\Bigr) d\alpha \Bigr] + \frac{T
h^2}{\pi^2 (4+h^2)} \int_{\beta}^{1}
\frac{\cos{(hL\alpha)}}{\alpha} d\alpha \\
+& \frac{8 T h}{\pi^2 (4+h^2)^2 L} \int_{\beta}^{1}
\frac{\sin{(hL\alpha)}}{\alpha^2} d\alpha + \frac{T}{\pi^2}
\int_{1}^{\infty} \frac{F(\alpha) - F_h(\alpha)}{\alpha^2} d\alpha
+ O(\beta^{-1} L^{-2}).
\end{split}
\end{equation*}
Apply (\ref{4.1}) and (\ref{4.2}) directly,
$$S_{3,2} = \frac{T}{8L^2 \beta^2} - \frac{T
\cos{(hL\beta)}}{2(4+h^2) L^2 \beta^2} + \frac{Th\sin{(hL\beta)}}
{(4+h^2)^2 L^3 \beta^3} + O(L^{-4}).$$ Similar to the treatment of
$S_{3,1}$,
$$S_{3,3} = \frac{T}{16 \pi^6 L^2 \beta^4} \int_{0}^{\infty}
(F(\alpha) - F_h(\alpha)) k''\Bigl(\frac{\alpha}{2\pi \beta}\Bigr)
d\alpha$$ and we split the integral into four pieces:
$$J = \int_{0}^{\infty} = \int_{0}^{\beta} +
\int_{\beta}^{1-\epsilon_T} + \int_{1-\epsilon_T}^{1} +
\int_{1}^{\infty} = J_1 + J_2 + J_3 + J_4.$$ By Lemma
\ref{lemma2.3}, \ref{lemma2.4} and \ref{lemma2.5},
\begin{equation*}
\begin{split}
J_1 =& 2(\pi^6 - 6\pi^4) \beta^2 \Bigl[1 - \frac{4\cos{(hL\beta)}}
{4+h^2}\Bigr] - \frac{16 \pi^4 h \sin{(hL\beta)}}{4+h^2} \beta^3 L
\\
&- \frac{32 \pi^4 h^2 \cos{(hL\beta)}}{(4+h^2)^2} \beta^2 +
\frac{16(\pi^6 - 4\pi^4) h \sin{(hL\beta)}}{(4+h^2)^2 L} \beta \\
&+ \frac{16 \pi^2 h^2}{4+h^2} \beta^2 L^2 \int_{0}^{\beta} \alpha
\cos{(hL\alpha)} k\Bigl(\frac{\alpha}{2\pi \beta}\Bigr) d\alpha
\\
&+ \frac{128 \pi^2 h}{(4+h^2)^2} \beta^2 L \int_{0}^{\beta}
\sin{(hL\alpha)} k\Bigl(\frac{\alpha}{2\pi \beta}\Bigr) d\alpha +
O(\beta^{-1} L^{-2}).
\end{split}
\end{equation*}
By Lemma \ref{lemma2.6} and \ref{lemma2.7},
\begin{equation*}
\begin{split}
J_2 + J_3 =& \int_{\beta}^{1} \Bigl[ \alpha - \frac{4
\cos{(hL\alpha)}}{4+h^2} \alpha + \frac{8h
\sin{(hL\alpha)}}{(4+h^2)^2 L} \Bigr] \frac{24\pi^4
\beta^4}{\alpha^4} d\alpha \\
=& 12 \pi^4 (\beta^2 - \beta^4) - \frac{16 \pi^4 (12+h^2)
\cos{(hL\beta)}}{(4+h^2)^2} \beta^2 \\
+& \frac{16 \pi^4 (12+h^2) \cos{(hL)}}{(4+h^2)^2} \beta^4 +
\frac{16h(12+h^2) \pi^4 \sin{(hL\beta)}}{(4+h^2)^2} \beta^3 L \\
-& \frac{16h(12+h^2) \pi^4 \sin{(hL)}}{(4+h^2)^2} \beta^4 L +
\frac{64 \pi^4 h \sin{(hL\beta)}}{(4+h^2)^2 L} \beta \\
-& \frac{64 \pi^4 h \sin{(hL)}}{(4+h^2)^2 L} \beta^4 + \frac{16
\pi^4 h^2 (12+h^2)}{(4+h^2)^2} \beta^4 L^2 \int_{\beta}^{1}
\frac{\cos{(hL\alpha)}}{\alpha} d\alpha + O(L^{-2}).
\end{split}
\end{equation*}
$$J_4 = 24 \pi^4 \beta^4 \int_{1}^{\infty} \frac{F(\alpha) -
F_h(\alpha)}{\alpha^4} d\alpha.$$ Thus,
\begin{equation*}
\begin{split}
S_{3,3} =& \frac{T}{8\beta^2 L^2} - \frac{3T}{4\pi^2 L^2} -
\frac{T\cos{(hL\beta)}}{2(4+h^2)L^2\beta^2} + \frac{8Th
\sin{(hL\beta)}}{\pi^2 (4+h^2)^2 \beta L} \\
&+ \frac{Th\sin{(hL\beta)}}{(4+h^2)^2 L^3 \beta^3} + \frac{T(12 +
h^2) \cos{(hL)}}{\pi^2 (4+h^2)^2 L^2} - \frac{Th(12+h^2)
\sin{(hL)}}{\pi^2 (4+h^2)^2 L} \\
& - \frac{4Th\sin{(hL)}}{\pi^2 (4+h^2)^2 L^3} + \frac{Th^2}{\pi^4
\beta^2 (4+h^2)} \int_{0}^{\beta} \alpha \cos{(hL\alpha)}
k\Bigl(\frac{\alpha}{2\pi \beta}\Bigr) d\alpha
\\
&+ \frac{8Th}{\pi^4 \beta^2 (4+h^2)^2 L} \int_{0}^{\beta}
\sin{(hL\alpha)} k\Bigl(\frac{\alpha}{2\pi \beta}\Bigr) d\alpha
\\
&+ \frac{T h^2(12+h^2)}{\pi^2 (4+h^2)^2} \int_{\beta}^{1}
\frac{\cos{(hL\alpha)}}{\alpha} d\alpha \\
&+ \frac{3T}{2\pi^2 L^2} \int_{1}^{\infty} \frac{F(\alpha) -
F_h(\alpha)}{\alpha^4} d\alpha + O(\beta^{-1} L^{-4}).
\end{split}
\end{equation*}
Finally, combining the results for $S_{3,1}$, $S_{3,2}$ and
$S_{3,3}$, we have
\begin{equation*}
\begin{split}
S_3 =& \frac{T}{\pi^2} \int_{0}^{\beta} \frac{1 -
\cos{(hL\alpha)}}{\alpha} \Bigl[1 - \frac{\pi \alpha}{2 \beta}
\cot{\bigl(\frac{\pi \alpha}{2 \beta} \bigr)}\Bigr]^2 d\alpha +
\frac{T}{\pi^2} \int_{\beta}^{1} \frac{1 -
\cos{(hL\alpha)}}{\alpha} d\alpha \\
+& \frac{T}{\pi^2} \int_{1}^{\infty} \frac{F(\alpha) -
F_h(\alpha)}{\alpha^2} d\alpha + \frac{Th\sin{(hL)}}{\pi^2 (4+h^2)
L} + \frac{3T}{4\pi^2 L^2} - \frac{T(12+h^2) \cos{(hL)}}{\pi^2
(4+h^2)^2 L^2} \\
-& \frac{3T}{2\pi^2 L^2} \int_{1}^{\infty} \frac{F(\alpha) -
F_h(\alpha)}{\alpha^4} d\alpha + O(\beta^{-1} L^{-2})
\end{split}
\end{equation*}
as
$$\int_{\beta}^{1} \frac{\sin{(hL\alpha)}}{\alpha^2} d\alpha =
-\sin{(hL)} + \frac{\sin{(hL\beta)}}{\beta} + hL \int_{\beta}^{1}
\frac{\cos{(hL\alpha)}}{\alpha} d\alpha.$$ Remark: We keep some of
the $O(TL^{-2})$ terms explicit because, with more effort, one can
make the error term $=C_1 T L^{-2} + o(T L^{-2})$.
\section{Proof of Theorem \ref{theorem1}}
Take $\beta = 1/2$. Combining the results on $S_3$, $S_4$ and
$S_5$, we have
\begin{equation*}
\begin{split}
& \int_{0}^{T} (S(t+h) - S(t))^2\, dt \\
=& \frac{T}{\pi^2} \int_{0}^{1} \frac{1 -
\cos{(hL\alpha)}}{\alpha} d\alpha + \frac{T}{\pi^2}
\int_{1}^{\infty} \frac{F(\alpha) - F_h(\alpha)}{\alpha^2} d\alpha
\\
+& \frac{T}{\pi^2} \Bigl[\log{\log{2}} + C_0 - \sum_{m=2}^{\infty}
\sum_{p} \frac{1}{mp^m} \Bigr] (1 - \cos{(h\log{2})}) -
\frac{T}{\pi^2} \int_{0}^{h\log{2}} \frac{1 -
\cos{\alpha}}{\alpha} d\alpha \\
-& \frac{Th}{\pi^2} \int_{2}^{\infty} r(u)
\frac{\sin{(h\log{u})}}{u} du + \frac{T}{\pi^2}
\sum_{m=2}^{\infty} \sum_{p} \frac{1 - \cos{(hm\log{p})}}{m^2 p^m}
\\
+& \frac{Th \sin{(hL)}}{\pi^2 (4+h^2)L} + \frac{T}{\pi^2 L^2}
\frac{h^2(20+3h^2)}{4(4+h^2)^2} -\frac{3T}{2\pi^2 L^2}
\int_{1}^{\infty} \frac{F(\alpha)-F_h(\alpha)}{\alpha^4} d\alpha +
O\Bigl(\frac{T}{L^2}\Bigr)
\end{split}
\end{equation*}
which gives the theorem. Again, one can make the error term $=C_1
T L^{-2} + o(T L^{-2})$ with more effort.

The theorem improves that of Fujii [\ref{F}] as
\begin{itemize}
\item $h$ is allowed to be $O(1)$. \item All the terms besides
first two contribute $O(T h^2) + O(T L^{-2})$. \item It is
conjectured in [\ref{C2}] that
\begin{equation}
\label{5.1} F_h(\alpha) = F(\alpha) \frac{4\cos{(hL
\alpha)}}{4+h^2} + o(1) \mbox{ for } 1 \leq \alpha \leq A \mbox {
with arbitrary large } A.
\end{equation}
\end{itemize}

\bigskip

Furthermore, as
\begin{equation*}
\begin{split}
& \sum_{p \leq x} \frac{1-\cos{(h\log{p})}}{p} = \int_{2}^{x}
\frac{1-\cos{(h\log{u})}}{u\log{u}}  du + \int_{2^-}^{x} 1 -
\cos{(h\log{u})} dr(u) \\
=& \int_{h\log{2}}^{h\log{x}} \frac{1-\cos{\alpha}}{\alpha}
d\alpha - r(2^-)(1-\cos{(h\log{2})}) \\
&- h \int_{2^-}^{\infty} r(u) \frac{\sin{(h\log{u})}}{u} du +
O\Bigl(\frac{\log{x}}{\sqrt{x}}\Bigr),
\end{split}
\end{equation*}
Theorem \ref{theorem1} gives
\begin{equation*}
\begin{split}
& \int_{0}^{T} (S(t+h) - S(t))^2\, dt \\
=& \frac{T}{\pi^2} \int_{0}^{1} \frac{1 -
\cos{(hL\alpha)}}{\alpha} d\alpha + \frac{T}{\pi^2}
\int_{1}^{\infty} \frac{F(\alpha) - F_h(\alpha)}{\alpha^2} d\alpha
\\
+& \frac{T}{\pi^2} \Bigl[\sum_{m=1}^{\infty} \sum_{p^m \leq x}
\frac{1 - \cos{(hm\log{p})}}{m^2 p^m} + Ci(h\log{x}) -
\log{(h\log{x})} - C_0\Bigr]
\\
+& \frac{Th \sin{(hL)}}{\pi^2 (4+h^2)L} + \frac{T}{\pi^2 L^2}
\frac{h^2(20+3h^2)}{4(4+h^2)^2} -\frac{3T}{2\pi^2 L^2}
\int_{1}^{\infty} \frac{F(\alpha)-F_h(\alpha)}{\alpha^4} d\alpha +
O\Bigl(\frac{T}{L^2}\Bigr)
\end{split}
\end{equation*}
where $Ci(x) = -\int_{x}^{\infty} \frac{\cos{t}}{t} dt = C_0 +
\log{x} + \int_{0}^{x} \frac{\cos{t}-1}{t} dt$ is the cosine
integral. If we assume Montgomery's conjecture [\ref{M}] on
$F(\alpha)$ and (\ref{5.1}), the first two terms account for the
GUE part of Berry's formula ($19$) conjectured in [\ref{B}] by a
similar calculation as page $79$ of [\ref{F}]. Moreover, the third
term is the non-GUE part of Berry's formula. So, our theorem is
even more precise than Berry's formula.
\section{Proof of Theorem \ref{theorem2}}
First, let us consider
$$\mbox{L}(x,t) = \sum_{0 < \gamma \leq T} \frac{x^{i (\gamma
- t)}}{1 + (t - \gamma)^2}$$ where the sum here is over of the
imaginary parts of the non-trivial zeros of the Riemann zeta
function.

\begin{lemma}
\label{lemma6.2} For all $\alpha$ and $h$, we have
$$F_h(\alpha) \leq F(\alpha).$$
\end{lemma}

Proof: First, by partial fractions and Cauchy's residue theorem,
$$\int_{-\infty}^{\infty} \frac{1}{(1+(t-a)^2)(1+(t-b)^2)} dt =
\frac{2\pi}{4 + (a-b)^2}.$$ Then
\begin{equation*}
\begin{split}
0 \leq& \int_{-\infty}^{\infty} |\mbox{L}(x,t) - \mbox{L}(x,t-h)|^2 dt \\
=& 2 \sum_{0 < \gamma, \gamma' \leq T} x^{i(\gamma - \gamma')}
\int_{-\infty}^{\infty} \frac{1}{(1 + (t - \gamma)^2)(1+(t -
\gamma')^2)} dt \\
&- 2Re \sum_{0 < \gamma, \gamma' \leq T} x^{i(\gamma - \gamma'-h)}
\int_{-\infty}^{\infty} \frac{1}{(1 + (t - \gamma)^2)(1+(t - h -
\gamma')^2)} dt \\
=& \pi F(x,T) - \pi F_h(x,T).
\end{split}
\end{equation*}
Set $x = (\frac{T}{2\pi})^\alpha$ and divide through by
$\frac{TL}{2\pi^2}$, we have the lemma.

\bigskip

Assuming {\bf RH}, Montgomery [\ref{M}] proved that, for fixed $0
< \beta < 1$,
\begin{equation}
\label{6.1} \Bigl(\frac{TL}{2\pi}\Bigr)^{-1} \sum_{0 < \gamma,
\gamma' \leq T} \biggl[ \frac{\sin{\frac{\beta
(\gamma-\gamma')L}{2}}}{\frac{\beta (\gamma-\gamma')L}{2}}
\biggr]^2 w(\gamma-\gamma') \sim \frac{1}{\beta} + \frac{\beta}{3}
\end{equation}
as $T \rightarrow \infty$. This also holds for $\beta = 1$ by
Goldston [\ref{G2}]. Using (\ref{2.4}), one can prove similarly
that for fixed $0 < \beta \leq 1$,
\begin{equation}
\label{6.2}
\begin{split}
& \Bigl(\frac{TL}{2\pi}\Bigr)^{-1} \sum_{0 < \gamma, \gamma' \leq
T} \biggl[ \frac{\sin{\frac{\beta
(\gamma-\gamma'-h)L}{2}}}{\frac{\beta
(\gamma-\gamma'-h)L}{2}} \biggr]^2 w(\gamma-\gamma'-h) \\
\sim& \frac{1}{\beta} + \frac{8\beta}{4+h^2}
\frac{2\frac{\sin{(hL\beta)}}{hL\beta} - 1 - \cos{(hL\beta)}}
{(hL\beta)^2}
\end{split}
\end{equation}
as $T \rightarrow \infty$ under {\bf RH} only (similar to
[\ref{G2}] or using Lemma $7$ of [\ref{GM}] in the argument of
Chan [\ref{C2}]). Note that $2\frac{\sin{x}}{x} - 1 - \cos{x} \leq
\frac{x^2}{6}$ by simply looking at their Taylor series. We also
need the following
\begin{lemma}
\label{lemma6.3} For any real number $c$,
\begin{equation*}
\begin{split}
& \int_{c-1}^{c+1} F_h(\alpha) (1 - |\alpha-c|) d\alpha \\
=& \Bigl(\frac{TL}{2\pi}\Bigr)^{-1} \sum_{0 < \gamma, \gamma' \leq
T} \cos{(L(\gamma - \gamma' -h)c)} \biggl[ \frac{\sin{\frac{
(\gamma-\gamma'-h)L}{2}}}{\frac{(\gamma-\gamma'-h)L}{2}}
\biggr]^2 w(\gamma-\gamma'-h), \\
& \int_{c-1}^{c+1} F(\alpha) (1 - |\alpha-c|) d\alpha \\
=& \Bigl(\frac{TL}{2\pi}\Bigr)^{-1} \sum_{0 < \gamma, \gamma' \leq
T} (\frac{T}{2\pi})^{i c (\gamma-\gamma')} \biggl[
\frac{\sin{\frac{(\gamma-\gamma')L}{2}}}{\frac{
(\gamma-\gamma')L}{2}} \biggr]^2 w(\gamma-\gamma').
\end{split}
\end{equation*}
\end{lemma}

Proof: The second one follows from the first one by setting $h=0$.
To prove the first one, we have, from the definition of
$F_h(\alpha)$,
\begin{equation*}
\begin{split}
& \int_{c-1}^{c+1} F_h(\alpha) (1 - |\alpha - c|) d\alpha =
\int_{-1}^{1} F_h(\alpha-c) (1-|\alpha|) d\alpha \\
=& \sum_{0 < \gamma, \gamma' \leq T} \int_{-1}^{1} \cos{(L(\gamma
- \gamma' -h)(\alpha -c))} (1-|\alpha|) d\alpha\, w(\gamma -
\gamma' - h) \\
=& \sum_{0 < \gamma, \gamma' \leq T} 2 \cos{(L(\gamma - \gamma'
-h)c)} \frac{1 - \cos{(L(\gamma - \gamma' -h))}}{((\gamma -
\gamma'- h)L)^2} w(\gamma-\gamma'-h)
\end{split}
\end{equation*}
by integration by parts. This gives the lemma as $\cos{2x} = 1 -
2\sin^2{x}$.

Note: (\ref{6.1}) and (\ref{6.2}) can be proved by setting $c=0$
in Lemma \ref{lemma6.3} and using the asymptotic formulas for
$F(\alpha)$ and $F_h(\alpha)$ like (\ref{2.4}).

\begin{lemma}
\label{lemma6.4} Assume {\bf RH}. For any $\epsilon > 0$ and $T$
sufficiently large,
$$0 \leq \int_{c}^{c+1} \bigl(F(\alpha) - F_h(\alpha)\bigr)\,
d\alpha \leq \frac{16}{3} + \epsilon$$ uniformly for any real
number $c$.
\end{lemma}

Proof: From Lemma \ref{lemma6.2}, we have $F(\alpha) - F_h(\alpha)
\geq 0$. This gives the lower bound as well as
$$\frac{1}{2} \int_{c-1/2}^{c+1/2} \bigl(F(\alpha) - F_h(\alpha)
\bigr)\, d\alpha \leq \int_{c-1}^{c+1} \bigl(F(\alpha) -
F_h(\alpha)\bigr) (1 - |\alpha -c|)\, d\alpha.$$ So, by Lemma
\ref{lemma6.3},
\begin{equation*}
\begin{split}
\int_{c-1/2}^{c+1/2} \bigl(F(\alpha) - F_h(\alpha)\bigr)\, d\alpha
\leq& 2 \Bigl(\frac{TL}{2\pi}\Bigr)^{-1} \sum_{0 < \gamma, \gamma'
\leq T} \biggl[ \frac{\sin{\frac{(\gamma-\gamma')L}{2}}}{\frac{
(\gamma-\gamma')L}{2}} \biggr]^2 w(\gamma-\gamma') \\
+& 2 \Bigl(\frac{TL}{2\pi}\Bigr)^{-1} \sum_{0 < \gamma, \gamma'
\leq T} \biggl[ \frac{\sin{\frac{
(\gamma-\gamma'-h)L}{2}}}{\frac{(\gamma-\gamma'-h)L}{2}} \biggr]^2
w(\gamma-\gamma'-h).
\end{split}
\end{equation*}
Now, using (\ref{6.1}) and (\ref{6.2}) with $\beta = 1$ and
$2\frac{\sin{x}}{x} - 1 - \cos{x} \leq \frac{x^2}{6}$, the right
hand side of the above inequality is
$$\leq 2\Bigl(\frac{4}{3} + \frac{\epsilon}{4}\Bigr) + 2\Bigl(
\frac{4}{3} + \frac{\epsilon}{4}\Bigr) = \frac{16}{3} + \epsilon$$
when $T$ is sufficiently large.

We are now in the position to prove Theorem \ref{theorem2}. By
Lemma \ref{lemma6.2} and \ref{lemma6.4},
$$0 \leq \int_{1}^{\infty} \frac{F(\alpha) - F_h(\alpha)}
{\alpha^2} d\alpha \leq \sum_{c=1}^{\infty} \frac{1}{c^2}
\int_{c}^{c+1} \bigl(F(\alpha) - F_h(\alpha)\bigr)\, d\alpha < 5.4
\frac{\pi^2}{6} < 9.$$ Similarly,
$$0 \leq \int_{1}^{\infty} \frac{F(\alpha) - F_h(\alpha)}
{\alpha^4} d\alpha \leq \sum_{c=1}^{\infty} \frac{1}{c^4}
\int_{c}^{c+1} \bigl(F(\alpha) - F_h(\alpha)\bigr)\, d\alpha < 5.4
\frac{\pi^4}{90} < 6.$$

Tsz Ho Chan\\
Case Western Reserve University\\
Mathematics Department, Yost Hall 220\\
10900 Euclid Avenue\\
Cleveland, OH 44106-7058\\
USA\\
txc50@po.cwru.edu

\end{document}